\documentclass[12pt,amssymb,amstex,draft]{article}

\usepackage{amsmath}
\usepackage{amssymb}
\usepackage{amsxtra}
\usepackage{latexsym}
\usepackage{ifthen}
\usepackage{amsthm}

\begin{document}

\centerline{\bf ON $L^1$-FUNCTIONS WITH A VERY SINGULAR BEHAVIOUR}

\vskip 15 pt \centerline{Alexander A. Kovalevsky}

\vskip 10 pt

\centerline{\footnotesize {\it Institute of Applied Mathematics
and Mechanics}}

\centerline{\footnotesize{\it Rosa Luxemburg St. 74, 83114
Donetsk, Ukraine}}

\vskip 25 pt

\baselineskip 25 pt

\noindent {\bf Abstract.} We give examples of $L^{1}$-functions
that are essentially unbounded on every non-empty open subset of
their domains of definition. We obtain such functions as limits of
weighted sums of functions with the unboundedly increasing number
of singular points lying at the nodes of standard compressible
periodic grids in $\Bbb R^n$. Moreover, we prove that the latter
(basic) functions possess properties of uniform integral
boundedness but do not have a pointwise majorant. Some
applications of the main results are given.

\vskip 5 pt


\vskip 5 pt


\vskip 30 pt

\noindent {\bf 1. Introduction}

In this article, we give examples of $L^{1}$-functions that are
essentially unbounded on every nonempty open subset of their
domains of definition. We obtain such functions as limits of
weighted sums of functions with the unboundedly increasing number
of singular points lying at the nodes of standard compressible
periodic grids in $\Bbb R^n$. Moreover, we prove that the latter
(basic) functions possess properties of uniform integral
boundedness but do not have a pointwise majorant.

The results obtained allow us to make some important conclusions
concerning the conditions under which the $\Gamma$-compactness of
integral functionals defined on variable weighted Sobolev spaces
was studied in [7,9,10]. However, we think that the main results
of the article are of a self-contained interest as well.

\vskip 30 pt

\noindent {\bf 2. Functions with singularities at the nodes of
periodic grids}

Let $n\in\mathbb{N}$, $n\geqslant 2$, and let $\Omega$ be a
bounded domain of $\mathbb{R}^n$.

For every $y \in \Bbb R^{n}$ and for every $\rho >0$ we set
$$
B(y,\rho)= \{x \in \Bbb R^{n}: |x-y|<\rho\},
$$
and for every $y \in \Bbb R^{n}$ and for every $t \in \Bbb N$ we
define
$$
Q_{t}(y)= \Big\{x \in \Bbb R^{n}: |x_{i}-y_{i}|< \frac{1}{2t},\,
\,i=1,\dots,n\Big\}.
$$
Moreover, for every $t \in \Bbb N$ we set
$$
Y_{t}=\{y \in \Bbb R^{n}: t y_{i} \in \Bbb Z, \,\,i=1,\dots,n\}.
$$

We have
$$
\forall \,t \in \Bbb N, \hskip 15 pt \bigcup_{y\in Y_t}\,
\overline{Q_t(y)} = \Bbb R^{n}, \eqno(2.1)
$$
$$
\forall\,t\in \Bbb N, \ \forall\,y,y^\prime \in Y_t,\,y \neq
y^\prime, \hskip 15 pt Q_{t}(y)\cap Q_{t}(y^\prime)=\emptyset.
\eqno(2.2)
$$
The proof of these assertions is simple.

Obviously, for every $t \in \Bbb N$ the elements of the set
$Y_{t}$ are the nodes of a periodic grid in $\Bbb R^{n}$ with the
period $1/t$. Such standard  grids are often used for instance in
different constructions to prove $\Gamma$-compactness of integral
functionals  and $G$-compactness of differential operators with
variable domain of definition (see for example [5,6]).

Next, for every $t \in \Bbb N$ we set
$$
X_{t}=\{y \in Y_{t} : Q_{t}(y) \subset \Omega\}.
$$
Owing to (2.1), there exists $m \in \Bbb N$ such that for every $t
\in \Bbb N$, $t>m$, we have $X_{t} \neq  \emptyset$.

For every $t \in \Bbb N$, $t>m$, we set
$$
G_{t}= \bigcup_{y \in X_t}\,B\Big(y, \frac{1}{2t}\Big).
$$
It is easy to see that for every $t \in \Bbb N$, $t>m$, and for
every $y \in X_{t}$, $B\big(y,\frac 1{2t}\big) \subset Q_t(y)
\subset \Omega$. Therefore, for every $t \in \Bbb N$, $t>m$, we
have $G_{t} \subset \Omega$.

Let for every $t \in \Bbb N$, $t>m$, $\chi_{t}:\Omega \to \Bbb R$
be the characteristic function of the set $G_{t}$, and let for
every $t \in \Bbb N$, $t>m$, and for every $y \in X_{t}$,
$\chi_{t,y}:\Omega \to \Bbb R$ be the characteristic function of
the ball $B\big(y,\frac{1}{2t}\big)$.

Further, we denote by $\mathcal K$ the set of all functions
$\sigma:[ 0, +\infty) \to (0, +\infty)$ with the properties:

\vskip 5 pt

\noindent \hskip 6.5 pt (i) $\,\sigma$ is continuous in $(0,
+\infty)$;

\vskip 5 pt

\noindent \hskip 3 pt (ii) $\,\sigma \geqslant 1$ in $[0,1]$;

\vskip 5 pt \noindent (iii) $\,\sigma(\rho)\to +\infty\,$ as
$\,\rho>0\,$ and $\,\rho \to 0$.

\vskip 10 pt

\noindent \hskip 1 pt (iv) $\,\displaystyle
\int_0^1\sigma(\rho)\rho^{n-1}d \rho < + \infty$.

\vskip 7 pt

For instance if $\,\sigma_{1}:[0,+\infty) \to (0,+\infty)$ is the
function such that $\sigma_{1}(0)=1$ and
\,$\sigma_{1}(\rho)=1/\rho$ for every $\rho>0$, then $\sigma_{1}
\in
 \mathcal{K}$.

For every $\sigma \in \mathcal{K}$ we set
$$
M_{\sigma}=\sigma(1)+n \int_0^1 \sigma(\rho) \rho^{n-1}d\rho.
$$

Let us give the following definition: if $\sigma \in \mathcal{K}$
and $t \in \Bbb N$, $t>m$, then $\nu^{(\sigma)}_{t}: \Omega \to
\Bbb R$ is the function such that for every $x \in \Omega$,
$$
\nu^{(\sigma)}_{t}(x)= \sigma(1)(1-\chi_{t}(x))+ \sum_{y \in
X_{t}} \sigma(2t|x-y|) \chi_{t,y}(x).
$$

\vskip 7 pt

{\bf Lemma 2.1.} {\it Let} $\,\sigma \in \mathcal{K}$ {\it and}
$\,t \in \Bbb N$, $t>m$. {\it Then the following properties hold}:

\vskip 5 pt

\noindent \hskip 6 pt (i) \,$\forall \,x \in \Omega \setminus
G_{t}$, \,\,$\nu^{(\sigma)}_{t}(x)= \sigma(1)$;

\vskip 5 pt

\noindent \hskip 3 pt (ii) \,{\it if}\, $\,y \in X_{t}$ {\it and}
$\,x \in B\big(y, \frac{1}{2t}\big)$, {\it then}
$\,\nu_{t}^{(\sigma)}(x)= \sigma(2 \,t |x-y|)$;

\vskip 5 pt

\noindent (iii) \,$\nu_{t}^{(\sigma)} \geqslant 1$ {\it in}
$\,\Omega$;

\vskip 5 pt

\noindent (iv) \,{\it the function} $\,\nu_{t}^{(\sigma)}$ {\it is
measurable};

\vskip 5 pt

\noindent \hskip 3 pt (v) \,$\nu_{t}^{(\sigma)} \in L^{1}(\Omega)$
{\it and} $\,\| \nu_{t}^{(\sigma)}\|_{L^{1}(\Omega)} \leqslant
M_{\sigma}\text{meas}\,\Omega$.

\vskip 5 pt

{\bf Proof.} Properties (i) and (ii) are immediate consequences of
the definition of the function $\nu^{(\sigma)}_{t}$. From these
properties, taking into account that $\sigma \geqslant 1$ in
$[0,1]$, we deduce property (iii). Moreover, using properties (i)
and (ii) and the continuity of $\sigma$ in $(0,+\infty)$, we
establish that the function $\nu^{(\sigma)}_{t}$ is continuous in
$\Omega \setminus X_{t}$. Therefore, the function
$\nu^{(\sigma)}_{t}$ is measurable. Thus property (iv) holds.

Next, by property (i), the function $\nu^{(\sigma)}_{t}$ is
summable in $\Omega \setminus G_{t}$ and
$$
\int_{\Omega \setminus G_{t}} \nu^{(\sigma)}_{t} dx = \sigma(1)
\text{meas}(\Omega \setminus G_{t}). \eqno(2.3)
$$
Moreover, taking into account property (ii) and the properties of
$\sigma$, by means of the change of variables, we establish that
for every $y \in X_{t}$ the function $\nu^{(\sigma)}_{t}$ is
summable in $B\big(y, \frac{1}{2t}\big)$ and
$$
\int_{B(y,\frac 1{2t})} \nu^{(\sigma)}_{t}dx =
\frac{\varkappa_{n}}{(2t)^{n}}\int_0^1 \sigma(\rho)
\rho^{n-1}d\rho, \eqno(2.4)
$$
where $\varkappa_{n}$ is the surface area of the unit sphere of
$\Bbb R^{n}$.

Now, taking into account (2.2), we conclude that the function
$\nu^{(\sigma)}_{t}$ is sum-mable in $\Omega$ and, by (2.3) and
(2.4),
\begin{multline}
\int_\Omega \nu^{(\sigma)}_{t}dx = \int_{\Omega \setminus G_{t}}
\nu^{(\sigma)}_{t}dx + \sum_{y \in X_{t}}\,\int_{B(y,
\frac 1{2t})} \nu_{t}^{(\sigma)}dx \notag \\[7pt]
= \sigma(1)\text{meas}(\Omega \setminus G_{t}) + |X_{t}|\,
\frac{\varkappa_{n}}{(2t)^{n}} \int_0^1 \sigma(\rho) \rho^{n-1}
d\rho, \tag{2.5}
\end{multline}
where $|X_{t}|$ is the number of elements of the set $X_{t}$.

From the definition of the set $G_{t}$ it follows that
$$
\text{meas}\,G_{t}= \frac{|X_{t}|}{(2t)^{n}}\,\text{meas}
\,B(0,1).
$$
This along with the equality $\varkappa_{n}=n\text{meas}\,B(0,1)$
and (2.5) implies that
$$
\int_\Omega \nu^{(\sigma)}_{t}dx = \sigma(1)\text{meas}(\Omega
\setminus G_{t})+ n\biggl( \int_0^1
\sigma(\rho)\rho^{n-1}d\rho\biggr)\text{meas}\,G_{t}.
$$
Hence we get the inequality
$\|\nu^{(\sigma)}_{t}\|_{L^{1}(\Omega)}\leqslant
M_{\sigma}\text{meas}\,\Omega$. Thus property (v) holds. \
$\square$

\vskip 7 pt

{\bf Remark 2.2.} If $\sigma \in \mathcal K$, $t \in \Bbb N$,
$t>m$, and $y \in X_{t}$, then $\nu^{(\sigma)}_{t}(x) \to +\infty$
as $x \in B\big(y, \frac{1}{2t}\big)\setminus \{y\}$ and $x \to
y$. This follows from property (ii) of Lemma 2.1 and the fact that
$\sigma(\rho) \to +\infty$ as $\rho>0$ and $\rho \to 0$.

\vskip 7 pt

Further, for every $\sigma \in \mathcal{K}$ and for every $t \in
\Bbb N$ we set
$$
\mu^{(\sigma)}_{t}= \sum_{k=1}^t\,k^{-2} \nu^{(\sigma)}_{m+k}\,.
$$

\vskip 7 pt

{\bf Lemma 2.3.} {\it Let} $\,\sigma \in \mathcal{K}$ {\it and}
$\,t \in \Bbb N$. {\it Then the following properties hold}:

\vskip 5 pt

\noindent (a) \,$\mu^{(\sigma)}_{t} \geqslant 1$ \ {\it in} \
$\Omega$;

\vskip 5 pt

\noindent (b) \,$\mu^{(\sigma)}_{t}< \mu^{(\sigma)}_{t+1}$ \ {\it
in} \ $\Omega$;

\vskip 5 pt

\noindent (c) \,$\mu^{(\sigma)}_{t} \in L^{1}(\Omega)$ \ {\it and}
\ $\|\mu^{(\sigma)}_{t}\|_{L^{1}(\Omega)} \leqslant
2M_{\sigma}\text{meas}\,\Omega$.

\vskip 5 pt

{\bf Proof.} Property (a) is a consequence of property (iii) of
Lemma 2.1. From the definitions of $\mu^{(\sigma)}_{t}$ and
$\mu^{(\sigma)}_{t+1}$ and property (iii) of Lemma 2.1 we deduce
property (b). Finally, property (c) follows from property (v) of
Lemma 2.1. \ $\square$

\vskip 30 pt

\noindent {\bf 3. Locally unbounded $L^1$-functions}

We denote by $\mathcal M$ the set of all functions $\mu\in
L^1(\Omega)$ with the properties:

\vskip 5 pt

\noindent \hskip 3 pt (i) \,$\mu\geqslant 1$ in $\Omega$;

\vskip 5 pt

\noindent (ii) \,for every nonempty open set $G \subset \Omega$
and for every $C>0$ there exists

\hskip 4 pt a measurable set $H\subset G$ such that
$\text{meas}\,H>0$ and $\mu\geqslant C$ in $H$.

\vskip 7 pt

{\bf Theorem 3.1.} {\it Let} $\,\sigma\in \mathcal{K}$. {\it Then
there exists a function} $\mu^{(\sigma)}\in \mathcal M$ {\it such
that}
\begin{align}
\mu^{(\sigma)}_t &\to \mu^{(\sigma)}\,\,\,\,\text{{\it
a.\,e. in}}\,\,\,\Omega, \tag{3.1} \\[5pt]
\|\mu^{(\sigma)}_t\|_{L^1(\Omega)}&\to
\|\mu^{(\sigma)}\|_{L^1(\Omega)}, \tag{3.2}\\[5pt]
\forall\,t\in \Bbb N,\hskip 15 pt \mu^{(\sigma)}_t &\leqslant
\mu^{(\sigma)}\,\,\,\,\text{{\it a.\,e. in}}\,\,\,\Omega.
\tag{3.3}
\end{align}

\vskip 5 pt

{\bf Proof.} By properties (b) and (c) of Lemma 2.3 and B.\,Levi's
theorem (see for instance [3, p.\,303]), there exists a function
$\tilde{\mu}\in L^1(\Omega)$ such that
\begin{align}
\mu^{(\sigma)}_t &\to \tilde{\mu}\,\,\,\,\text{a.\,e. in}\,\,\,
\Omega, \tag{3.4} \\[7pt]
\int_\Omega \mu^{(\sigma)}_t dx &\to \int_\Omega \tilde{\mu}\,dx.
\tag{3.5}
\end{align}
According to (3.4), there exists a set $E\subset \Omega$ with
measure zero such that
$$
\forall\,x\in \Omega\setminus E, \hskip 15 pt \mu^{(\sigma)}_t(x)
\to \tilde{\mu}(x). \eqno(3.6)
$$

We define the function $\mu^{(\sigma)}:\Omega\rightarrow \Bbb R$
by
$$
\mu^{(\sigma)}(x)=
\begin{cases}
\tilde{\mu}(x) & \text{if}\quad x\in \Omega\setminus E,\\
1 &\text{if}\quad x\in E.
\end{cases}
$$
Clearly, $\mu^{(\sigma)}\in L^1(\Omega)$. Using (3.6) and property
(a) of Lemma 2.3, we establish that $\mu^{(\sigma)}\geqslant 1$ in
$\Omega$. Moreover, owing to (3.6), assertion (3.1) holds, and due
to (3.5), assertion (3.2) is valid. Besides, by (3.6) and property
(b) of Lemma 2.3, assertion (3.3) holds.

Next, let $G$ be a nonempty open set of $\Bbb R^n$, $G\subset
\Omega$, and let $C>0$. We fix $z\in G$. Obviously, there exists
$\rho_0>0$ such that
$$
B(z,\rho_0)\subset G. \eqno(3.7)
$$
We fix $l\in \Bbb N$ such that $l>n/\rho_0$. By (2.1), there
exists $y\in Y_{m+l}$ such that $z\in \overline{Q_{m+l}(y)}$.
Since $n/l<\rho_0$, we have $Q_{m+l}(y)\subset B(z,\rho_0)$. This
and (3.7) imply that
$$
Q_{m+l}(y)\subset G. \eqno(3.8)
$$
Therefore, $Q_{m+l}(y)\subset \Omega$. Hence
$$
y\in X_{m+l}. \eqno(3.9)
$$

Since $\sigma\in \mathcal K$, we have $\sigma(\rho)\to +\infty$ as
$\rho>0$ and $\rho\to 0$. Then there exists $\delta\in (0,1)$ such
that
$$
\forall\,\rho\in (0,\delta), \hskip 15 pt \sigma(\rho)>(C+1)l^2.
\eqno(3.10)
$$
We set
$$
G^\prime=B\Big(y,\frac{\delta}{2(m+l)}\Big)\setminus \{y\}.
$$
Evidently, $G^\prime\subset Q_{m+l}(y)$. From this and (3.8) we
get $G^\prime\subset G$.

Let $x\in G^\prime$. We have $2(m+l)|x-y|\in (0,\delta)$.
Therefore, by (3.10),
$$
\sigma(2(m+l)|x-y|)>(C+1)l^2. \eqno(3.11)
$$
Moreover, taking into account (3.9) and using property (ii) of
Lemma 2.1, we obtain
$$
\nu^{(\sigma)}_{m+l}(x)=\sigma(2(m+l)|x-y|). \eqno(3.12)
$$
Finally, by the definition of $\mu^{(\sigma)}_l$ and property
(iii) of Lemma 2.1, we have
$$
\mu^{(\sigma)}_l(x)\geqslant l^{-2} \nu^{(\sigma)}_{m+l}(x).
\eqno(3.13)
$$

From (3.11)--(3.13) we infer that
$$
\forall\,x\in G^\prime, \hskip 15 pt \mu_l^{(\sigma)}(x) > C+1.
\eqno(3.14)
$$

Further, by assertion (3.1) and D.\,Egoroff's theorem (see for
instance [3, p.\,287]), there exists a measurable set
$\Omega^\prime\subset \Omega$ such that
$$
\text{meas}(\Omega\setminus \Omega^\prime)\leqslant \frac 12
\,\text{meas}\,G^\prime, \eqno(3.15)
$$
$$
\mu^{(\sigma)}_t \to \mu^{(\sigma)}\,\,\,\,\text{uniformly in
}\,\,\Omega^\prime. \eqno(3.16)
$$
We set $H=G^\prime\cap\Omega^\prime$. Clearly, the set $H$ is
measurable and $H\subset G$. Moreover, $G^\prime\subset
H\cup(\Omega\setminus \Omega^\prime)$. This and (3.15) imply that
$\,\text{meas}\,G^\prime\leqslant \text{meas}\,H + \frac
12\,\text{meas}\,G^\prime$. Hence $\text{meas}\,H>0$.

According to (3.16), there exists $t_0\in \Bbb N$ such that for
every $t\in \Bbb N$, $t\geqslant t_0$, and for every $x\in
\Omega^\prime$,
$$
|\mu^{(\sigma)}_t(x)-\mu^{(\sigma)}(x)|\leqslant 1. \eqno(3.17)
$$
Let $x\in H$. We fix $t\in \Bbb N$, $t\geqslant \max(l,t_0)$.
Using (3.17), property (b) of Lemma 2.3 and (3.14), we obtain
$\mu^{(\sigma)}(x)\geqslant \mu^{(\sigma)}_t(x)-1\geqslant
\mu^{(\sigma)}_l(x)-1>C$. Thus $\mu^{(\sigma)}\geqslant C$ in $H$.
Now, we conclude that $\,\mu^{(\sigma)}\in \mathcal M$. \
$\square$

\vskip 7 pt

By virtue of Theorem 3.1, the set $\mathcal{M}$ is nonempty. Let
us state several propositions describing some properties of this
set.

\vskip 7 pt

{\bf Proposition 3.2.} {\it For every} $\lambda>1$ {\it we have}
$\mathcal{M}\cap L^\lambda(\Omega)\neq \emptyset$.

\vskip 5 pt

{\bf Proof.} Let $\lambda>1$. We fix $\mu\in \mathcal M$ and set
$\mu_\lambda=\mu^{1/\lambda}$. Since $\mu\in \mathcal M$, we have
$\mu_\lambda\in \mathcal M$. Moreover, due to the definition of
$\mu_\lambda$ and the inclusion $\mu\in L^1(\Omega)$, we have
$\mu_\lambda\in L^\lambda(\Omega)$. Thus $\mu_\lambda\in \mathcal
M\cap L^\lambda(\Omega)$. Hence $\mathcal M\cap
L^\lambda(\Omega)\neq \emptyset$. \ $\square$

\vskip 7 pt

{\bf Proposition 3.3.} {\it Let} $F:(0,+\infty)\to (0,+\infty)$
{\it be a nondecreasing continuous function, and} $F(1)=1$. {\it
Let} $\,\sigma\in \mathcal K$. {\it Suppose that} $\,\sigma
F(\sigma) \not\in \mathcal K$. {\it Then there exists a function}
$\mu\in \mathcal M$ {\it such that} $\mu F(\mu)\not\in
L^1(\Omega)$.

\vskip 5 pt

{\bf Proof.} We set $\sigma_\ast=\sigma F(\sigma)$. Due to the
properties of $F$ and the inclusion $\sigma\in \mathcal K$, the
function $\sigma_\ast$ has the following properties: $\sigma_\ast$
is continuous in $(0,+\infty)$, $\,\sigma_\ast\geqslant 1$ in
$[0,1]$ and $\sigma_\ast(\rho)\to +\infty$ as $\rho>0$ and
$\rho\to 0$. Hence, taking into account that $\sigma_\ast\not\in
\mathcal K$, we obtain that
$$
\int^1_0 \sigma_\ast(\rho) \rho^{n-1}d\rho=+\infty. \eqno(3.18)
$$

Next, by Theorem 3.1, there exists a function $\mu\in \mathcal M$
such that
$$
\mu^{(\sigma)}_1\leqslant \mu \,\,\,\,\text{a.\,e. in }\,\,\Omega.
\eqno(3.19)
$$

Suppose that
$$
\mu F(\mu)\in L^1(\Omega). \eqno(3.20)
$$
Using (3.19), the definition of $\mu^{(\sigma)}_1$ and the fact
that $F$ is nondecreasing, we establish that
$\,\nu^{(\sigma)}_{m+1} F(\nu^{(\sigma)}_{m+1})\leqslant \mu
F(\mu)\,$ a.\,e in $\,\Omega$. This and (3.20) imply that
$\nu^{(\sigma)}_{m+1} F(\nu^{(\sigma)}_{m+1})\in L^1(\Omega)$.

Now, we fix $y\in X_{m+1}$ and for every $\varepsilon\in (0,1)$
set
$$
K_\varepsilon=\Big\{x\in \Bbb R^n:\frac{\varepsilon}{2(m+1)}<
|x-y|<\frac{1}{2(m+1)}\Big\}.
$$
Obviously, for every $\varepsilon\in (0,1)$,
$$
\int_{K_\varepsilon} \nu^{(\sigma)}_{m+1}
F(\nu^{(\sigma)}_{m+1})dx\leqslant \|\nu^{(\sigma)}_{m+1}
F(\nu^{(\sigma)}_{m+1})\|_{L^1(\Omega)}. \eqno(3.21)
$$
On the other hand, using property (ii)  of Lemma 2.1, the
definition of $\sigma_\ast$ and the change of variables, we obtain
that for every $\varepsilon\in (0,1)$,
$$
\int_{K_\varepsilon} \nu^{(\sigma)}_{m+1}
F(\nu^{(\sigma)}_{m+1})dx = \frac{\varkappa_n}{2^n(m+1)^n}
\int^1_\varepsilon \sigma_\ast(\rho)\rho^{n-1}d\rho.
$$
This and (3.18) imply that
$$
\int_{K_\varepsilon} \nu^{(\sigma)}_{m+1}
F(\nu^{(\sigma)}_{m+1})dx\to +\infty\,\,\,\,\text{as} \,\,\,\,
\varepsilon\to 0.
$$
However, the result obtained contradicts (3.21). Due to this
contradiction, we conclude that inclusion (3.20) does not hold.
Thus $\mu F(\mu)\not\in L^1(\Omega)$. \ $\square$

\vskip 7 pt

{\bf Corollary 3.4.} {\it Let} $\lambda>0$. {\it Then there exists
a function} $\mu\in \mathcal M$ {\it such that} $\mu(\ln
\mu)^\lambda\not\in L^1(\Omega)$.

\vskip 5 pt

{\bf Proof.} Let $F:(0,+\infty)\to (0,+\infty)$ be the function
such that for every $\rho\in (0,+\infty)$, $\,
F(\rho)=[\,\ln(e-1+\rho)]^\lambda$. Clearly, the function $F$ is
nondecreasing and continuous, and $F(1)=1$.

Let $\sigma:[0,+\infty)\to (0,+\infty)$ be the function such that
$$
\sigma(\rho)=
\begin{cases}
\displaystyle \frac{4}{\rho^n}\Big(\ln \frac{1}{\rho}\Big)^{-1}
\Big(\ln\ln
\frac{1}{\rho}\Big)^{-2} & \text{if} \quad 0<\rho<e^{-e},\\[10pt]
4e^{n e-1}& \text{if} \quad \rho=0 \,\,\,\,\text{or}\,\,\,\,
\rho\geqslant e^{-e}.
\end{cases}
$$
It is easy to see that the function $\sigma$ is continuous in
$(0,+\infty)$. In addition, we have
$$
\forall\,\rho\in(0,e^{-e}), \hskip 15 pt \sigma(\rho)\geqslant
\frac{1}{4\rho}\,. \eqno(3.22)
$$
Using this fact, we establish that $\sigma>1$ in $[0,1]$ and
$\sigma(\rho)\to +\infty$ as $\rho>0$ and $\rho\to 0$. Finally, it
is not difficult to verify that
$$
\int_0^1 \sigma(\rho)\rho^{n-1}d\rho< +\infty.
$$
The described properties of the function $\sigma$ allow us to
conclude that $\sigma\in \mathcal K$.

Now, let us show that $\sigma F(\sigma)\not\in \mathcal K$. In
fact, let $\rho\in(0,e^{-e})$. By (3.22) and the definition of
$F$, we have
$$
F(\sigma(\rho))>\Big(\ln \frac{1}{4\rho}\Big)^\lambda. \eqno(3.23)
$$
Since
$$
\ln \frac{1}{4\rho}= -2\ln 2+ \frac{2}{e}\ln \frac{1}{\rho}+
\Big(1- \frac{2}{e}\,\Big)\ln \frac{1}{\rho}
$$
and $\ln(1/\rho)>e$, the following inequality holds:
$$
\ln \frac{1}{4\rho}>\Big(1-\frac{2}{e}\Big)\ln \frac{1}{\rho}\,.
\eqno(3.24)
$$
Moreover, we observe that
$$
\Big(\ln\ln \frac{1}{\rho}\Big)^2<\frac{4}{\lambda^2}\Big(\ln
\frac{1}{\rho}\Big)^\lambda. \eqno(3.25)
$$

From (3.23)--(3.25) we deduce that for every $\rho\in (0,e^{-e})$,
$$
\sigma(\rho)F(\sigma(\rho))
\rho^{n-1}>\frac{\lambda^{2}}{\rho}\Big(1-\frac{2}{e}\,\Big)^\lambda
\Big(\ln \frac{1}{\rho}\Big)^{-1}.
$$
Hence
$$
\int^1_0 \sigma(\rho) F(\sigma(\rho))\rho^{n-1}d\rho= +\infty.
$$
Therefore, $\sigma F(\sigma)\not\in \mathcal{K}$. Then, by
Proposition 3.3, there exists a function $\mu\in \mathcal{M}$ such
that $\mu F(\mu)\not\in L^1(\Omega)$. Hence, taking into account
that $\mu\in L^1(\Omega)$ and $\,\mu F(\mu)\leqslant
2^\lambda\mu+2^\lambda\mu(\ln \mu)^\lambda \,\,\,\,\text{in}\,\,\,
\Omega\,$, we infer that $\mu(\ln\mu)^\lambda\not\in L^1(\Omega)$.
\ $\square$

\vskip 7 pt

{\bf Corollary 3.5.} {\it There exists a function} $\mu\in
\mathcal{M}$ {\it such that for every} $\lambda>1$, $\mu\not\in
L^\lambda(\Omega)$.

\vskip 5 pt

{\bf Proof.} By Corollary 3.4, there exists a function $\mu\in
\mathcal{M}$ such that
$$
\mu \ln\mu \not\in L^1(\Omega). \eqno(3.26)
$$
Let $\lambda>1$. Since $(\lambda-1)\ln \mu<\mu^{\lambda-1}$ in
$\Omega$, we have $\,\mu\ln \mu< \frac{1}{\lambda-1}\,\mu^\lambda
\,\,\,\,\text{in}\,\,\,\Omega$. This and (3.26) imply that
$\mu\not\in L^\lambda(\Omega)$. \ $\square$

\vskip 30 pt

\noindent {\bf 4. An exhaustion property of the domain $\Omega$}

In this section, we establish that unions of certain balls
connected with all the sets $X_t$, $t>m$, exhaust the domain
$\Omega$. This property is essentially used in Section 5 to study
the pointwise behaviour of the functions $\nu_t^{(\sigma)}$.

We set $\,\alpha=2^{-n}\text{meas}\,B(0,1)$. Evidently, $\alpha\in
(0,1)$.

For every $k$, $t\in \Bbb N$ we set
$$
\mathcal{B}^{(k)}_t=\bigcup_{y\in X_{m+t}}
B\Big(y,\frac{1}{2(m+t)k}\Big).
$$
Clearly, if $k, t\in \Bbb N$, we have $\mathcal{B}^{(k)}_t\subset
\Omega$.

\vskip 5 pt

{\bf Proposition 4.1.} {\it Let} $\,k\in \Bbb N$. {\it Then for
every open set} $G\subset \Omega$ {\it we have}
$$
\liminf_{t\to \infty}\,\text{meas}(G\cap
\mathcal{B}_t^{(k)})\geqslant\alpha k^{-n}\text{meas}\,G.
\eqno(4.1)
$$

\vskip 5 pt

{\bf Proof.} Let $G$ be an open set of $\Bbb R^n$ such that
$G\subset \Omega$. In the case $G=\emptyset$ inequality (4.1) is
obvious. Consider the case $G\neq \emptyset$. We fix
$\varepsilon>0$ and for every $j\in \Bbb N\,$ set $G_j=\{x\in G:
d(x,\partial\,G)> 1/j\}$. Clearly, $\text{meas}\,G_j \to
\text{meas}\,G$. Therefore, there exists $l\in \Bbb N$ such that
$G_l\neq \emptyset$ and
$$
\text{meas}(G\setminus G_l)\leqslant \varepsilon. \eqno(4.2)
$$
We fix $t\in \Bbb N$ such that $t\geqslant nl$ and set
$$
X^\prime_t=\{y\in Y_{m+t}:Q_{m+t}(y)\cap G_l \neq \emptyset\}.
$$
Moreover, we denote by $q_t$ the number of elements of the set
$X^\prime_t$. By (2.1), we have $X^\prime_t\neq \emptyset$, and
owing to (2.1) and the inequality $n/t\leqslant 1/l$, we get
$$
G_l\subset \bigcup_{y\in X^\prime_t} \overline{Q_{m+t}(y)}\subset
G. \eqno(4.3)
$$
This and (2.2) imply that
$$
(m+t)^{-n} q_t\leqslant \text{meas}\,G. \eqno(4.4)
$$

Next, from the obvious inclusion $G\setminus
\mathcal{B}^{(k)}_t\subset (G_l\setminus \mathcal{B}^{(k)}_t)\cup
(G\setminus G_l)$ and (4.2) we obtain
$$
\text{meas}(G\setminus \mathcal{B}^{(k)}_t)\leqslant
\text{meas}(G_l\setminus \mathcal{B}^{(k)}_t)+\varepsilon.
\eqno(4.5)
$$
Let us estimate the measure of the set $G_l\setminus \mathcal
B^{(k)}_t$. First of all we observe that, due to (4.3) and the
inclusion $G\subset \Omega$,
$$
X^\prime_t\subset X_{m+t}. \eqno(4.6)
$$
Let $x\in G_l\setminus \mathcal{B}_t^{(k)}$. Since $x\in G_l$, by
(4.3), there exists $y\in X^\prime_t$ such that $x\in
\overline{Q_{m+t}(y)}$. At the same time $x\not\in B\big(y,
\frac{1}{2(m+t)k}\big)$. This follows from (4.6) and the fact that
$x\not\in \mathcal{B}^{(k)}_t$. Thus $x\in
\overline{Q_{m+t}(y)}\setminus B\big(y,\frac{1}{2(m+t)k}\big)$,
and we conclude that
$$
G_l\setminus \mathcal{B}^{(k)}_t\subset \bigcup_{y\in X^\prime_t}
\Big[\,\overline{Q_{m+t}(y)}\setminus B\Big(y,
\frac{1}{2(m+t)k}\Big)\Big].
$$
Hence
$$
\text{meas}(G_l\setminus \mathcal{B}^{(k)}_t)\leqslant (1-\alpha
k^{-n})(m+t)^{-n} q_t. \eqno(4.7)
$$
From (4.5), (4.7) and (4.4) we deduce that \ $\text{meas}
(G\setminus \mathcal{B}^{(k)}_t)\leqslant (1-\alpha k^{-n})
\text{meas}\,G$ $+ \varepsilon$. Therefore, $\text{meas}(G\cap
\mathcal{B}^{(k)}_t)\geqslant \alpha k^{-n}\text{meas}\, G -
\varepsilon$. Hence we get (4.1). \ $\square$

\vskip 7 pt

{\bf Corollary 4.2.} {\it Let} $\,k \in \Bbb N$. {\it Then for
every measurable set} $H\subset \Omega$ {\it we have}
$$
\liminf_{t\to \infty}\,\text{meas}(H \cap
\mathcal{B}^{(k)}_t)\geqslant \alpha k^{-n}\text{meas}\,H.
\eqno(4.8)
$$

\vskip 5 pt

{\bf Proof.} Let $H$ be a measurable set of $\Bbb R^n$ such that
$H\subset \Omega$. We fix $\varepsilon>0$. Clearly, there exists
an open set $H^\prime$ of $\Bbb R^n$ such that $H^\prime\subset
\Omega$ and
$$
\text{meas}(H\setminus H^\prime)< \varepsilon, \hskip 15 pt
\text{meas}(H^\prime\setminus H)< \varepsilon. \eqno(4.9)
$$
By Proposition 4.1, we have
$$
\liminf_{t\to \infty}\,\text{meas}(H^\prime\cap
\mathcal{B}^{(k)}_t)\geqslant \alpha k^{-n}\text{meas}\,H^\prime.
$$
This and (4.9) imply that
$$
\liminf_{t\to \infty}\,\text{meas}(H\cap
\mathcal{B}^{(k)}_t)\geqslant \alpha k^{-n}(\text{meas}\,H -
\varepsilon)- \varepsilon.
$$
Hence we get (4.8). \ $\square$

The following result describes the above-mentioned exhaustion
property of the domain $\Omega$.

\vskip 5 pt

{\bf Proposition 4.3.} {\it For every} $k\in \Bbb N$ {\it we have}
$$
\text{meas}\Big(\Omega\setminus \bigcup^\infty_{t=1}\,
\mathcal{B}^{(k)}_t\Big)=0. \eqno(4.10)
$$

\vskip 5 pt

{\bf Proof.} Let $k\in \Bbb N$. We set
$$
\Phi=\Omega\setminus \bigcup_{t=1}^\infty\,\mathcal{B}^{(k)}_t.
$$
By Corollary 4.2, we have
$$
\alpha k^{-n}\text{meas}\,\Phi\leqslant \liminf_{t\to
\infty}\,\text{meas}(\Phi \cap \mathcal{B}^{(k)}_t), \eqno(4.11)
$$
and from the definition of $\Phi$ it follows that for every $t\in
\Bbb N$, $\,\Phi\cap \mathcal{B}^{(k)}_t=\emptyset$. The latter
fact and (4.11) imply that $\text{meas}\,\Phi=0$. Thus equality
(4.10) holds. \ $\square$

\vskip 7 pt

{\bf Remark 4.4.} The exhaustion property described by Proposition
4.3 is an analogue of the exhaustion condition which is assumed in
some results of [4,8].

\vskip 30 pt

\noindent {\bf 5. Further properties of the functions
$\nu^{(\sigma)}_t$}

{\bf Theorem 5.1.} {\it Let} $\,\sigma\in \mathcal{K}$. {\it Then
for almost every} $x\in \Omega$ {\it the sequence}
$\{\nu^{(\sigma)}_{m+t}(x)\}$ {\it is unbounded}.

\vskip 5 pt

{\bf Proof.} For every $k\in \Bbb N$ we set
$$
\mathcal{B}^{(k)}=\bigcup_{t=1}^\infty\,\mathcal{B}^{(k)}_t.
$$
Then we define
$$
E_0= \bigcup_{k=1}^\infty\,(\Omega\setminus \mathcal{B}^{(k)}).
$$
From Proposition 4.3 it follows that $\text{meas}\,E_0=0$.

Next, we set
$$
E_1=\bigcup_{t=1}^\infty\,X_{m+t}.
$$
Clearly, $\text{meas}\,E_1=0$.

We fix $x\in \Omega\setminus (E_0\cup E_1)$. Suppose that
$$
\text{the sequence}\,\, \{\nu^{(\sigma)}_{m+t}(x)\}\,\,\text{is
 bounded}. \eqno(5.1)
$$
Then there exists $M>0$ such that
$$
\forall\,t\in \Bbb N, \hskip 15 pt
\nu^{(\sigma)}_{m+t}(x)\leqslant M. \eqno(5.2)
$$
Since $\sigma\in \mathcal{K}$, we have $\sigma(\rho)\to +\infty$
as $\rho>0$ and $\rho\to 0$. Therefore, there exists $\delta>0$
such that
$$
\forall\,\rho\in (0,\delta), \hskip 15 pt \sigma(\rho)>M.
\eqno(5.3)
$$
We fix $j\in \Bbb N$ such that $j>1/\delta$. Since $x\in
\Omega\setminus E_0$, we have $x\in \mathcal{B}^{(j)}$. Then there
exists $l\in \Bbb N$ such that $x\in \mathcal{B}^{(j)}_l$. Hence,
taking into account the definition of $\mathcal{B}^{(j)}_l$, we
obtain that there exists $y\in X_{m+l}$ such that $x\in B\big(y,
\frac{1}{2(m+l)j}\big)$. Therefore, by property (ii) of Lemma 2.1,
we have
$$
\nu^{(\sigma)}_{m+l}(x)= \sigma(2(m+l)|x-y|). \eqno(5.4)
$$
Moreover, $2(m+l)|x-y|<1/j<\delta$. At the same time, due to the
fact that $x\not\in E_1$, we have $x\neq y$. Thus $2(m+l)|x-y|\in
(0,\delta)$. This along with (5.3) and (5.4) implies that
$\nu^{(\sigma)}_{m+l}(x)>M$. However, by (5.2), we have
$\nu^{(\sigma)}_{m+l}(x)\leqslant M$. The contradiction obtained
proves that assertion (5.1) is not valid. Therefore, the sequence
$\{\nu^{(\sigma)}_{m+t}(x)\}$ is unbounded. \ $\square$

\vskip 7 pt

{\bf Corollary 5.2.} {\it Let} $\,\sigma\in \mathcal{K}$. {\it
Then there is no function} $\,\psi:\Omega\to \Bbb R$ {\it such
that for every} $\,t\in \Bbb N$, $\,\nu^{(\sigma)}_{m+t}\leqslant
\psi\,$ {\it a.\,e. in} $\,\Omega$.

\vskip 5 pt

{\bf Proof.} Suppose that there exists a function $\psi: \Omega\to
\Bbb R$ such that for every $t\in \Bbb N$,
$\,\nu^{(\sigma)}_{m+t}\leqslant \psi\,$ a.\,e. in $\,\Omega$.
Then there exists a set $E^\prime\subset \Omega$ with measure zero
such that
$$
\text{for every}\,\,\, x\in \Omega\setminus E^\prime
\,\,\,\text{and for every}\,\,\, t\in \Bbb N \,\,\,\text{we
have}\,\,\, \nu^{(\sigma)}_{m+t}(x)\leqslant \psi(x). \eqno(5.5)
$$
Moreover, by Theorem 5.1, there exists a set
$E^{\prime\prime}\subset \Omega$ with measure zero such that
$$
\text{for every}\,\,\, x\in \Omega\setminus E^{\prime\prime}
\,\,\,\text{the sequence}\,\,\,
\{\nu^{(\sigma)}_{m+t}(x)\}\,\,\,\text{is unbounded}. \eqno(5.6)
$$
Let $x\in \Omega\setminus (E^\prime\cup E^{\prime\prime})$. Then,
in view of (5.5), the sequence $\{\nu^{(\sigma)}_{m+t}(x)\}$ is
bounded. At the same time, by (5.6), the sequence
$\{\nu^{(\sigma)}_{m+t}(x)\}$ is unbounded. The contradiction
obtained leads to the conclusion required. \ $\square$

\vskip 7 pt

{\bf Theorem 5.3.} {\it Let} $\,\sigma\in \mathcal{K}$. {\it Then
for every open cube} $Q\subset \Bbb R^n$ {\it we have}
$$
\limsup_{t\to \infty}\,\int_{Q\cap \Omega}
\nu^{(\sigma)}_{m+t}\,dx \leqslant M_\sigma\text{meas}(Q \cap
\Omega). \eqno(5.7)
$$

\vskip 5 pt

{\bf Proof.} Let $Q$ be an open cube of $\Bbb R^n$. If $Q \cap
\Omega=\emptyset$, inequality (5.7) is evident. Thus we may
consider that $Q\cap \Omega \neq \emptyset$.

We have
$$
Q=\{x\in \Bbb R^n:|x_i-z_i|<a/2, \,\,\, i=1,\dots,n\},
$$
where $z\in \Bbb R^n$ and $a>0$.

We fix $\varepsilon\in (0,1)$ and set
$$
Q_\varepsilon=\{x\in \Bbb R^n: |x_i-z_i|<(1+\varepsilon)a/2,
\,\,\, i=1,\dots,n\}.
$$
It is easy to see that $Q\subset Q_\varepsilon$ and
$$
\text{meas}(Q_\varepsilon\setminus Q)\leqslant (2a)^n
n\varepsilon. \eqno(5.8)
$$

Next, we fix $z^\prime\in Q\cap \Omega$. Clearly, there exists
$\rho_0>0$ such that
$$
B(z^\prime,\rho_0)\subset Q\cap \Omega. \eqno(5.9)
$$
We fix $t\in \Bbb N$ such that
$t>\max\big\{\frac{n}{\rho_0},\frac{2}{a\varepsilon}\big\}$ and
set
$$
\tilde{X}_t=\{y\in X_{m+t}:Q \cap Q_{m+t}(y)\neq \emptyset\}.
$$
Observe that $\tilde{X}_t\neq \emptyset$. In fact, by (2.1), there
exists $y\in Y_{m+t}$ such that $z^\prime\in
\overline{Q_{m+t}(y)}$. Evidently, $Q \cap Q_{m+t}(y)\neq
\emptyset$. Moreover, if $x\in Q_{m+t}(y)$, we have
$|x-z^\prime|\leqslant |x-y|+ |z^\prime-y|<n/t<\rho_0$. This and
(5.9) imply that $Q_{m+t}(y)\subset \Omega$. Hence $y\in X_{m+t}$.
Now, we may conclude that $y\in \tilde{X}_t$. Therefore, the set
$\tilde{X}_t$ is nonempty.

We denote by $\tilde{q}_t$ the number of elements of the set
$\tilde{X}_t$. Since
$$
\bigcup_{y\in \tilde{X}_t}Q_{m+t}(y)\subset Q_\varepsilon \cap
\Omega,
$$
using (2.2) and (5.8), we get
$$
(m+t)^{-n}\,\tilde{q}_t \leqslant \text{meas}(Q\cap \Omega)+
(2a)^n n\varepsilon. \eqno(5.10)
$$

Further, we set
$$
G^\prime_t=(Q\cap \Omega)\setminus G_{m+t}, \hskip 15 pt
G^{\prime\prime}_t=\bigcup_{y\in \tilde{X}_t} B\Big(y,
\frac{1}{2(m+t)}\Big).
$$
It is easy to see that $Q\cap \Omega\subset G^\prime_t \cup
G^{\prime\prime}_t$, $\,G^\prime_t\subset \Omega$ and
$G^{\prime\prime}_t\subset \Omega$. Then
$$
\int_{Q\cap \Omega} \nu^{(\sigma)}_{m+t}\,dx \leqslant
\int_{G^\prime_t}\nu^{(\sigma)}_{m+t}\,dx +
\int_{G^{\prime\prime}_t}\nu^{(\sigma)}_{m+t}\,dx. \eqno(5.11)
$$
Taking into account property (i) of Lemma 2.1, we get
$$
\int_{G^\prime_t}\nu^{(\sigma)}_{m+t}\,dx \leqslant
\sigma(1)\text{meas}(Q\cap \Omega), \eqno(5.12)
$$
and using (2.4) and (5.10), we obtain
\begin{multline}
\int_{G^{\prime\prime}_t}\nu^{(\sigma)}_{m+t}\,dx= \sum_{y\in
\tilde{X}_t}\,\int_{B(y,
\frac{1}{2(m+t)})}\nu^{(\sigma)}_{m+t}\,dx=\frac{\varkappa_n
\,\tilde{q}_t}{2^n(m+t)^n}\int_0^1 \sigma(\rho)\rho^{n-1}d\rho
\notag \\[7pt]
\leqslant \bigg(n\int_0^1
\sigma(\rho)\rho^{n-1}d\rho\bigg)[\,\text{meas}(Q\cap \Omega)+
(2a)^n n\varepsilon\,]. \tag{5.13}
\end{multline}
From (5.11)--(5.13) it follows that
$$
\int_{Q\cap \Omega} \nu^{(\sigma)}_{m+t}\,dx \leqslant M_\sigma
[\,\text{meas}(Q\cap \Omega)+ (2a)^n n\varepsilon\,].
$$
Hence we deduce (5.7). \ $\square$

\vskip 30 pt

\noindent {\bf 6. Some applications}

Let $p\in (1,n)$. We denote by $\mathcal{N}_p$ the set of all
nonnegative functions $\nu:\Omega\to \Bbb R$ such that $\nu>0$
a.\,e. in $\Omega$, $\nu\in L^1_\text{loc}(\Omega)$ and
$(1/\nu)^{1/(p-1)}\in L^1_{\text{loc}}(\Omega)$.

Observe that $\mathcal{M}\subset \mathcal{N}_p$.

In [7,9,10] some weighted Sobolev spaces $W_s$ associated with the
exponent $p$, a weight $\nu\in \mathcal{N}_p$ and a sequence of
domains $\Omega_s\subset \Omega$ were considered, and theorems on
the $\Gamma$-compactness of the sequence of integral functionals
$J_s:W_s\to \Bbb R$ of the form
$$
J_s(u)=\int_{\Omega_s}f_s(x,\nabla u)dx
$$
were established.

Here we are not giving the corresponding definitions and
statements of the results of the above-mentioned articles. We only
point to several things connected with the conditions under which
the $\Gamma$-compactness of integral functionals $J_s$ was proved.

In the above-mentioned articles, it is supposed that the
integrands $f_s:\Omega_s\times \Bbb R^n\to \Bbb R$ of the
functionals $J_s$ satisfy the following conditions:

\noindent (a$_1$) for every $s\in \Bbb N$ and for every $\xi\in
\Bbb R^n$ the function $f_s(\cdot,\xi)$ is measurable

\hskip 7 pt in $\Omega_s$;

\noindent (a$_2$)  for every $s\in \Bbb N$ and for almost every
$x\in \Omega_s$ the function $f_s(x,\cdot)$ is

\hskip 7 pt convex in $\Bbb R^n$;

\noindent (a$_3$) for every $s\in \Bbb N$, for almost every $x\in
\Omega_s$ and for every $\xi\in \Bbb R^n$,

\hskip 7 pt $c_1\nu(x)|\xi|^p - \psi_s(x)\leqslant f_s(x,\xi)
\leqslant c_2\nu(x)|\xi|^p + \psi_s(x)$.

\vskip 5 pt In the latter condition $c_1$ and $c_2$ are positive
constants and $\{\psi_s\}$ is a sequence of functions such that

\noindent (b$_1$) for every $s\in \Bbb N$, $\psi_s\in
L^1(\Omega_s)$ and $\psi_s \geqslant 0$ in $\Omega_s$;

\noindent (b$_2$) for every open cube $Q\subset \Bbb R^n$,
$$
\limsup_{s\to \infty}\,\int_{Q\cap \Omega_s} \psi_s\,dx \leqslant
\int_{Q\cap \Omega} b\,dx,
$$
where $b\in L^1(\Omega)$ and $b\geqslant 0$ in $\Omega$.

According to results of Section 5, there exist sequences that
satisfy conditions (b$_1$) and (b$_2$) but do not have pointwise
majorants. Indeed, the following simple proposition holds.

\vskip 5 pt

{\bf Proposition 6.1.} {\it Let} $\,\sigma\in \mathcal{K}$. {\it
Let} $\,b:\Omega\to \Bbb R$ {\it be the function such that for
every} $x\in \Omega$, $\,b(x)=M_\sigma$. {\it Let for every} $s\in
\Bbb N$, $\psi_s=\nu^{(\sigma)}_{m+s}$ {\it and}
$\Omega_s=\Omega$. {\it Then the sequence} $\{\psi_s\}$ {\it
satisfies conditions} (b$_1$) {\it and} (b$_2$) {\it but there is
no function} $\,\psi:\Omega\to \Bbb R$ {\it such that for every}
$s\in \Bbb N$, $\,\psi_s\leqslant \psi\,$ {\it a.\,e. in}
$\,\Omega$.

This result follows from properties (iii) and (v) of Lemma 2.1,
Theorem 5.3 and Corollary 5.2.

We note that Proposition 6.1 is of interest to compare conditions
(a$_1$)--(a$_3$) with conditions under which the
$\Gamma$-compactness of sequences of integral functionals with
degenerate variable integrands and the same domain of integration
was established in [1,2].

In [1] it is supposed that the integrands $g_s:\Bbb R^n\times \Bbb
R^n \to \Bbb R$ of the functionals under consideration satisfy
conditions of measurability and convexity like (a$_1$) and (a$_2$)
and the following condition:
\begin{align}
&\text{for every}\,\, s\in \Bbb N,\,\,\text{for almost every}\,\,
x\in
\Bbb R^n \,\,\text{and for every}\,\, \xi\in \Bbb R^n, \notag \\
&w_s(x)|\xi|^p \leqslant g_s(x,\xi)\leqslant \Lambda
w_s(x)(1+|\xi|^p), \tag{6.1}
\end{align}
where $\Lambda>0$ and $\{w_s\}$ is a sequence of nonnegative
functions on $\Bbb R^n$ satisfying a uniform Muckenhoupt
condition.

In order to compare condition (a$_3$) with condition (6.1), we
give the following example.

\vskip 5 pt

{\bf Example 6.2.} Suppose that $p\geqslant 2$ and all the
conditions of Proposition 6.1 are satisfied. Let for every $s\in
\Bbb N$ the function $f_s:\Omega \times \Bbb R^n\to \Bbb R$ be
defined by
$$
f_s(x,\xi)=\nu(x)|\xi|^p +
(\nu(x))^{(p-1)/p}(\psi_s(x))^{1/p}|\xi|^{p-1}, \hskip 12 pt
(x,\xi)\in \Omega \times \Bbb R^n. \eqno(6.2)
$$
It is easy to see that the sequence $\{f_s\}$ satisfies conditions
(a$_1$)--(a$_3$). At the same time, by Proposition 6.1, the
sequence $\{\psi_s\}$ satisfies conditions (b$_1$) and (b$_2$).

Assume that there exist $\lambda>0$ and a sequence of functions
$\varphi_s:\Omega\to \Bbb R$ such that
\begin{align}
&\text{for every}\,\,s\in \Bbb N, \,\, \text{for almost every}\,\,
x\in
\Omega \,\, \text{and for every}\,\, \xi\in \Bbb R^n, \notag \\
&\varphi_s(x)|\xi|^p\leqslant f_s(x,\xi)\leqslant \lambda
\varphi_s(x)(1+|\xi|^p). \tag{6.3}
\end{align}
This and the property $\nu>0\,$ a.\,e. in $\Omega$ imply that
there exists a set $\tilde{E}\subset \Omega$ with measure zero
such that
\begin{align}
&\forall\,x\in \Omega\setminus \tilde{E}, \hskip 15 pt \nu(x)>0,
\tag{6.4} \\
\text{for every}\,\, &s\in \Bbb N, \,\,\,\text{for every}\,\,x\in
\Omega\setminus
\tilde{E}\,\, \text{and for every} \,\, \xi\in \Bbb R^n, \notag \\
&\varphi_s(x)|\xi|^p\leqslant f_s(x,\xi)\leqslant \lambda
\varphi_s(x)(1+|\xi|^p). \tag{6.5}
\end{align}

We fix $s\in \Bbb N$ and $x\in \Omega\setminus \tilde{E}$. Let
$\xi\in \Bbb R^n$, $\xi\neq 0$. Using (6.5) and (6.2), we obtain
$$
\varphi_s(x)|\xi|^p\leqslant \nu(x)|\xi|^p +
(\nu(x))^{(p-1)/p}(\psi_s(x))^{1/p}|\xi|^{p-1} \leqslant
2\nu(x)|\xi|^p + \psi_s(x).
$$
Therefore, $\varphi_s(x)\leqslant 2\nu(x) + \psi_s(x)|\xi|^{-p}$.
Hence, passing to the limit as $|\xi|\to \infty$, we get
$$
\varphi_s(x)\leqslant 2\nu(x). \eqno(6.6)
$$
Now, let $\xi\in \Bbb R^n$, $|\xi|=1$. Using (6.2), (6.5) and
(6.6), we obtain
$$\nu(x)|\xi|^p +
(\nu(x))^{(p-1)/p}(\psi_s(x))^{1/p}|\xi|^{p-1}\leqslant \lambda
\varphi_s(x)(1+|\xi|^p) \leqslant 2\lambda \nu(x)(1+|\xi|^p).
$$
Therefore, $(\nu(x))^{(p-1)/p}(\psi_s(x))^{1/p}\leqslant 4\lambda
\nu(x)$. Hence, taking into account (6.4), we get
$\,\psi_s(x)\leqslant (4\lambda)^p \nu(x)$.

Thus for every $s\in \Bbb N$, $\,\psi_s\leqslant (4\lambda)^p
\nu\,$ \,a.\,e. in $\,\Omega$. However, this contradicts the fact
that, by Proposition 6.1, there is no function $\psi:\Omega\to
\Bbb R$ such that for every $s\in \Bbb N$, $\,\psi_s\leqslant
\psi\,$ a.\,e. in $\,\Omega$. The contradiction obtained proves
that there is no $\lambda>0$ and sequence of functions
$\varphi_s:\Omega\to \Bbb R$ such that assertion (6.3) holds.

As a result, we conclude that the sequence $\{f_s\}$ satisfies
conditions (a$_1$)--(a$_3$) but any extensions $g_s$ of the
functions $f_s$ on $\Bbb R^n \times \Bbb R^n$ do not satisfy
condition (6.1). Consequently, the sequence $\{f_s\}$ cannot be
considered in the framework of conditions imposed on the
integrands of functionals in [1]. The same conclusion concerns
conditions imposed on the integrands of functionals in [2].

Finally, mention should be made of the following. In [7,9,10] the
$\Gamma$-compact-ness of integral functionals  was proved under
the assumption that there exists a sequence of nonempty open sets
$\Omega^{(k)}$ of $\Bbb R^n$ such that for every $k\in \Bbb N\,$,
$\,\overline{\Omega^{(k)}}\subset \Omega^{(k+1)}\subset \Omega$,
$\text{meas}(\Omega\setminus\Omega^{(k)})\to 0$ and for every
$k\in \Bbb N$ the functions $\nu$ and $b$ are bounded in
$\Omega^{(k)}$. Evidently, if $\nu\in \mathcal{M}$, the given
assumption cannot be realized.

\vskip 100 pt

E-mail addresses:

alexkvl@iamm.ac.donetsk.ua

 aakovalevsky@yahoo.com

\vskip 15 pt Postal address:

Alexander A. Kovalevsky

Institute of Applied Mathematics and Mechanics

Rosa Luxemburg St. 74

83114 Donetsk

Ukraine

\vskip 15 pt

tel.: +38 062 311 04 53

fax:  +38 062 311 02 85

\end{document}